\documentclass[journal,onecolumn,draftclsnofoot,]{IEEEtran}
\usepackage{amssymb, amsmath, bm, graphicx, color, colortbl, xfrac, tabularx}
\usepackage{todonotes}
\DeclareGraphicsExtensions{.pdf,.png,.jpg, .eps}

\begin{document}

\title{Families of delta-correlated arrays derived from the Fourier spectra of Huffman sequences}

\author{Imants~D.~Svalbe*,
        David~M.~Paganin,
        Timothy~C.~Petersen%
  \thanks{All authors are with the School of Physics and Astronomy, Monash University, Victoria 3800, Australia.}
  \thanks{Timothy~C.~Petersen is primarily with the Monash Centre for Electron Microscopy, Monash University, Victoria 3800, Australia.}
}


\maketitle

\begin{abstract}
Finite discrete Huffman sequences, together with their extension to n-dimensional arrays, are highly valued because their discrete aperiodic auto-correlations optimally approximate the continuum form of the delta function. We present here several new families of real and integer-valued Huffman sequences, beyond those of the recursive form found by Hunt and Ackroyd. These new sequences are derived using the remarkably uniform discrete Fourier power spectra of Huffman sequences, where the elements are expressed as terms of the Fibonacci sequence.
\end{abstract}

\maketitle
\IEEEpeerreviewmaketitle

\section{Introduction}

A discrete (Kronecker) delta function $\delta_j$, comprising a single non-zero entry for element $j$ among an array of zeros, has practical significance in defining the limit of element localization, acting as does a Dirac delta $\delta$ on the continuum. Functions $f$ and $g$ where $f \cdot g = \delta$ under some binary operation ($\cdot$) are widely useful, as then $f$ and $g$ are functionally inverse.  For the correlation operator (where $\cdot \rightarrow \otimes$) in the periodic domain, it is possible to find functions $f$ that are self-inverse, such that $f$ represents a factorization of the delta function. Discrete periodic examples that approximate this behavior are $M$-, Legendre and non-redundant arrays, while Singer sequences perfectly correlate to a Kronecker delta \cite{Singer1938,McFarland2011}. Importantly, a sequence $f=f_1,f_2,\cdots,f_N$ of $N$ discrete elements with perfect delta auto-correlation in the periodic domain \cite{Luke1988}, where the right extended element $f_{N+1} = f_{1}$ and the left extended element $f_{0} = f_{N}$, will, in general, exhibit significantly lower auto-correlation under zero-padded \emph{aperiodic} conditions, where $f_{0} = f_{N+1} = 0$.

In the discrete aperiodic domain, the `end' auto-correlation terms are the product between the first element, $f_1$ and last element $f_N$. For any sequence with $N > 1$, this product will always be non-zero, even though the sum of the products at all other non-zero shifts may be zero. Hence any delta formed under aperiodic conditions must exhibit end entries in the sea of zeros around the central spike that denotes the limit of localization. Huffman \cite{Huffman1962} recognized and formalized this critical difference between periodic and aperiodic deltas. Throughout, we refer to Huffman sequences which have no more than one pair of non-zero off-peak aperiodic auto-correlation elements as `canonical' or `delta-correlated'. Such arrays meet the limiting bounds for off-peak correlation values established by Welch \cite{Welch1974} and Sarwate \cite{Sarwate1979,Sarwate1999}.

In addition to results from direct methods, Hunt and Ackroyd \cite{HuntAckroyd1980} found a recursive formula to produce desired integer-valued canonical Huffman sequences where the deltas formed under aperiodic auto-correlation have end entries with a minimal value of -1 (where $f_1 f_N = -1$ for $f$ of length $N$). Their approach placed the zeros of z-transforms on a single pair of circles with requisite radii $X$ and $1/X$, at the necessary regularly-spaced phase angles in the Argand plane.  Multiplication by carefully placed poles was used to create particular configurations of zeros to fully satisfy Huffman's criteria for canonical sequences \cite{Huffman1962}. A second order recursion, intrinsic to these arrays, was identified upon inspection from a series-expansion of the z-transform, evident also in their direct method solution and consistent with the coefficients of a transformation involving ultra-spherical polynomials \cite{HuntAckroyd1980}. With their scaling parameter set to unity, the z-transform $X$ and $1/X$ radii are interestingly equal to the golden ratio and its reciprocal.  This is perhaps unsurprising, given the intrinsic recursion matches that which defines Fibonacci polynomials.  As such, we anticipated that an elementary proof of the canonical delta-correlation might be possible using well known recursion identities.  

Our recent work \cite{svalbe2019sharp} showed that the single recursion identified by Hunt and Ackroyd can indeed express such canonical Huffman arrays as largely comprised of Fibonacci series that scale according to Fibonacci polynomials, for a continuous scaling parameter $s$. In that work, families of sequences that have delta-like correlation were constructed as discrete projections from higher dimensional ($nD$) canonical Huffman sequences, these are referred to as `quasi-Huffman arrays'  \cite{svalbe2019sharp}.  By connecting delta-correlation with the elementary diffraction physics of certain aberrated waves on the continuum, the flatness of the Fourier spectrum was identified as an important and useful property of Huffman sequences, which deserves further analysis \cite{svalbe2019sharp}.

Huffman \cite{Huffman1962}, Ackroyd \cite{Ackroyd1972,Ackroyd1975}, Schroeder \cite{Schroeder1970} and Golay \cite{Golay1972,Golay1975b} constructed similar $1D$ sequences of length $N$ using different Fourier techniques to approximate canonical delta-correlation, without precisely adhering to Huffman's criteria. They selected $N$ points in the $2D$ Fourier plane, distributed around the perimeter of two circles, with radii $X$ and $1/X$, with $X$ close to $1$. In these construction algorithms, the value chosen for $X$, the division of $N$ points being either inner or outer radii points, and the angular placement of these `root’ points around the perimeter, are sequence design variables. Given appropriate choices for these variables (often found by iteration after applying symmetry constraints), the resulting sequences can closely reproduce the canonical Huffman auto-correlation. A similar approach is mirrored in more recent work \cite{OjedaHuffman1994}. This method is advantageous as, by construction and with the value $X \approx 1$, it produces sequence examples where all element amplitudes are significantly more uniform than for the Fibonacci based Huffman sequences. However, truly uniform (binary) Huffman sequences do not exist \cite{White1977}. Most of the sequences generated as above are comprised of complex, non-integer elements (cf. Barker sequences \cite{BorweinJedwab2004} generalized to the complex plane \cite{Popovic1992, ChangGolombBarkerGener1996, FriesePolyphBarker1996}). 

The canonical Fibonacci-based Huffman sequences of Hunt and Ackroyd \cite{HuntAckroyd1980} are, by design, inherently real and integer valued for integer scaling parameters. While integer forms were ruled out in that work for lengths $N=4n+1$, one can find asymptotically close approximations for large $N=4n+1$ or exact canonical Gaussian-integer arrays, though we shall not consider complex arrays further here.  Also, as briefly given in the Supplementary Material Sec.~\ref{sec:nonFibonacci}, there exist length $4n+3$ canonical integer-valued arrays that are not based on the Fibonacci sequence, yet have identical correlation metrics, and direct methods can readily find $N = 4n+1$ integer-valued canonical Huffman arrays with end elements equal to 1.  The Supplementary Material Sec.~\ref{sec:Huffman proof} gives a formal and elementary proof that the aperiodic auto-correlation of the Fibonacci-based Huffman sequences is a discrete delta function, by summing auto-correlation products using standard recursion identities. 

Here we use key spectral properties of the aperiodic auto-correlation of Huffman sequences to show that the explicit Fourier transform of these Fibonacci-based Huffman sequences can be expressed in closed-form, using elementary methods.  The simplicity of these spectral insights prompts the discovery of many more canonical Huffman sequences that are not based upon the Fibonacci sequence, including new integer-valued arrays, as detailed in the latter part of this paper. 

The results presented here are for $1D$ sequences. In \cite{svalbe2019sharp}, the aforementioned extensions of $1D$ to $2D$ (see also \cite{Golomb1982}) and $nD$ Huffman arrays were presented, as were constructions for `quasi-Huffman' arrays, where the correlations closely approximate the exact $1D$ aperiodic form. The $1D$ proofs and results shown here can be readily extended for $nD$ Huffman arrays.

\section{The delta auto-correlation of Fibonacci based Canonical Huffman arrays} \label{sec:deltaAuto}
Using elementary means and recursion identities, it is proven in the Supplementary Material (Sec.~\ref{sec:Huffman proof}) that Huffman arrays comprising ordered sequences of Fibonacci polynomials $F_1(s), F_2(s), \cdots$ with the following form are canonically delta-correlated:
\begin{align}\label{CanonAutoCorr}
&H^N(s) = 2 s[(2s)^{-1},F_{1}(s),F_{2}(s),\cdots, F_{M}(s),F_{M+1}(s)/2 - F_{M}(s)/s,F_{-M}, \cdots,F_{-2}(s),F_{-1}(s),-(2 s)^{-1}], 
\end{align}
where $M$ must be even. The aperiodic auto-correlation $A_d(s)$ of Eq.~\ref{CanonAutoCorr} is zero for all off-peak shifts $|d| < N-1$, and non-zero at the end points $A_{N-1}(s) = A_{-(N-1)}(s) = -1$.  

The Fibonacci-based canonical Huffman family in Eq.~\ref{CanonAutoCorr} is identical to that derived by Hunt and Ackroyd and, with different choices of $s$, can be matched with all of their sequences found by direct methods \cite{HuntAckroyd1980}.

We begin our analysis of these interesting canonical Huffman sequences by demonstrating a closed-form expression for the auto-correlation peak at zero-shift, $A_0(s)$, using an elementary approach.  To this end, the central element of any auto-correlation is the sum of squared elements (by definition and in accord with Parseval's theorem), so the Huffman auto-correlation peak $A_0(s)$ can be written as, 
%
\begin{equation}\label{AutoZero}
A_0(s) = 2+ (2s)^2[(F_{M+1}(s)/2 - F_{M}(s)/s)^2
+ 2 \sum_{i=1}^{M}F_i(s)^2], \end{equation}
%
with $M = (N-3)/2$.  Equation~\ref{AutoZero} can be simplified using a generalized Fibonacci sum of squares identity \cite{Monk2004}, such that
\begin{align}\label{AutoZeroII}
&A_0(s) = 2+(2s)^2[(F_{M+1}(s)/2 - F_{M}(s)/s)^2 + 2 F_M(s)F_{M+1}(s)/s] \notag \\ 
&= 2+(2s)^2[(F_{M+1}(s)/2)^2 + (F_{M}(s)/s)^2 + F_M(s)F_{M+1}(s)/s] = 2+s^2F_{M+1}^2+2^2F_M F_{M+2},
\end{align}
where the Fibonacci polynomial recursion was used to arrive at the final equality.
We can now succinctly express the delta auto-correlation array $A(s)$ for canonical Huffman arrays $H^N(s)$:
\begin{equation}\label{HuffAutoCorrFinal}
A(s) = H^{N} \otimes H^{N}= [-1,0,\cdots,0,A_0,0,\cdots,0,-1],
\end{equation}
where $A_0=2+s^2F_{M+1}^2+2^2F_M F_{M+2}$.

For \emph{periodic} correlation, it is possible to devise `perfect arrays' with delta-auto-correlation (i.e.~all off-peak elements zero) \cite{Luke1988}.  For similar sequences that approximate perfect arrays, it is 
useful to compute quality metrics that compare to the delta-function \cite{Moharir1975}. Useful quality metrics include the relative sum of square auto-correlation elements as a `merit factor' $\mathcal{M}$, together with the peak to maximum off-peak ratio $\mathcal{R}$ and the degree of spectral-flatness $\mathcal{S}$.  These metrics are respectively defined as
\begin{equation}
    \mathcal{M} = \frac{A_0^2}{\sum_{i \ne 0} A_i^2},
\end{equation}
\begin{equation}
    \mathcal{R} = A_0/ \max(|A_i|_{i \ne 0}),
\end{equation}
\begin{equation}
    \mathcal{S}  =  \Delta F/\langle F \rangle,
\end{equation}
where $\Delta F$ is the maximum variation in the Fourier transform magnitudes, over all frequencies in the Fourier spectrum, and $\langle F \rangle$ is the mean Fourier magnitude.

Given the results detailed in this section, the first of these two quality metrics can be expressed in closed-form for canonical Huffman arrays:
\begin{equation}
    \mathcal{M} = (2+s^2F_{M+1}^2+2^2F_M F_{M+2})^2/2,
\end{equation}
\begin{equation}
    \mathcal{R} = 2+s^2F_{M+1}^2+2^2F_M F_{M+2}.
\end{equation}
The numeric values for metrics $\mathcal{M}$ and $\mathcal{R}$ are large compared to unity for canonical Huffman arrays, even for small $N$, $s = 1$.  The spectral flatness metric $\mathcal{S}$ is covered in the next section.

\section{The flat Fourier spectra of canonically delta-correlated arrays}

The Fourier transform $\mathcal{F}$ of a delta function is a constant, so that the power-spectrum is broad-band and flat.  The sparse delta-like form of the canonical Huffman auto-correlation means that the Fourier spectrum of a Huffman sequence $H$ has a similarly compact expression, which can be deduced using the Fourier convolution theorem applied to the auto-correlation, 
\begin{equation}\label{ConvTheorem}
\mathcal{F}[H \otimes H] = |\mathcal{F}[H]|^2.
\end{equation}
This section specifies the exact form of the power-spectrum for canonical Huffman arrays.  The simple delta-like form of the canonical delta auto-correlation of $H^N(s)$ foreshadows that its power-spectrum should be almost constant.  

Using a convention consistent with Eq.~\ref{ConvTheorem}, the discrete Fourier transform $g$ for sequence $f$ is here defined as,
\begin{equation}\label{DiscreteFourier}
g_q=\mathcal{F}[f]_q = \sum_{n=0}^{N-1}f_n e^{-2 \pi i n q/N},
\end{equation}
with inverse,
\begin{equation}\label{DiscreteInverseFourier}
f_n=\mathcal{F}^{-1}[g]_n = N^{-1} \sum_{q=0}^{N-1}g_q e^{2 \pi i q n/N}.
\end{equation}
A sequence of length $N = 2M+3$ has a bandwidth of frequencies ranging from zero up to the Nyquist frequency $M+1$. Upon taking the inverse Fourier transform of Eq.~\ref{ConvTheorem}, the \emph{periodic} auto-correlation $A^p(s)$ of $H^N(s)$ can be defined. In that case, the $-1$ values at the ends of the length $2N-1$ aperiodic auto-correlation in Eq.~\ref{HuffAutoCorrFinal} fold/wrap such that
\begin{equation}\label{HuffAutoCorrPeriodic}
A^p(s) = \mathcal{F}^{-1}[|\mathcal{F}[H^{N}(s)]|^2] = [A_0,-1,0,\cdots,0,-1].
\end{equation}
The $N$ wrapped terms in Eq.~\ref{HuffAutoCorrFinal} offer a potential short-cut to our proof that all off-peak correlations, bar the end points, are zero.  However this approach hides the remarkable summation to zero of cross products for partially overlapped sequences.

Before determining the power-spectrum of Eq.~\ref{HuffAutoCorrPeriodic} in closed-form, we shall evaluate the Fourier transform of the Huffman array at $q = 0$. As Eq.~\ref{DiscreteFourier} implies, the Fourier amplitude $\mathcal{F}[H^N(s)]_0$ at the zero-frequency $q = 0$ corresponds to the sum of the Huffman sequence, in which the $\pm 1$ end points and all even elements cancel.  By adding all of the odd elements of $H(s)$, for scale $s = 1$ (when $F_j(s) = f_j(s)$) we have $\underline x = -f_{M-2}$ and $\mathcal{F}[H^N(s)]_0 = 4 f_M - f_{M-2}$, for $s=1$, where $\underline x$ is the middle term in the sequence.  This relation generalizes for all $s$, in terms of Fibonacci polynomials as,
\begin{align}\label{Zeroq}
\mathcal{F}[H^N(s)]_0 &= \sum_{i=1}^{N} H^N_i(s) = s F_{M+1}(s) - 2 F_{M}(s) + 4 s \sum_{i=1}^{M/2} F_{2i-1}(s) \\ \nonumber
& = s F_{M+1}(s) + 2 F_{M}(s) = F_{M+2}(s) + F_{M}(s).
\end{align}
Above,  we have used the odd-index identity $\sum_{n=1}^N F_{2n-1}(s) = F_{2N}(s)/s$.  This identity is practically the same as that from the Fibonacci series ($s=1$), arising from the Fibonacci polynomial recursion $F_{2n-1}(s) = [F_{2n}(s)-F_{2(n-1)}(s)]/s$, such that the sum over $n$ is telescoping (consecutive terms cancel, leaving $F_{2N}(s)/s-F_0(s)/s$, with $F_0(s)=0$). 

We now return to evaluating the power-spectrum of canonical Huffman arrays, in the periodic domain, using Eq.~\ref{ConvTheorem}.  The Fourier transform of the periodic auto-correlation is given by substituting $A^p(s)$ in Eq.~\ref{HuffAutoCorrPeriodic} into Eq.~\ref{DiscreteFourier},
\begin{align}\label{AutoFourier}
\mathcal{F}[A^p(s)]_q = \sum_{n=0}^{N-1}A^{p}_n(s)
e^{-2 \pi i n q/N} = A_0 - 2 \cos (2 \pi q/N)
\end{align}
where $A_0 =2+s^2F_{M+1}^2+2^2F_M F_{M+2}$ is the peak of the auto-correlation $A^p(s)$.

Before concluding this section, we draw attention to an interesting observation that the $q=0$ origin in Eq.~\ref{AutoFourier}, $A_0-2$, appears different from the square of the zeroth order Fourier coefficient (in the last line of Eq.~\ref{Zeroq}). To reconcile this apparent contradiction, consider the square of Eq.~\ref{Zeroq} with the following recursion substitutions,  
\begin{align}\label{ReconcileSumSquare}
&|\mathcal{F}[H^N(s)]|^2_0 = [F_{M+2}(s) + F_{M}(s)]^2 = [F_{M+2}(s) + (F_{M+2}(s)-s F_{M+1}(s))]^2 = [2F_{M+2}(s) - s F_{M+1}(s)]^2 \notag \\
&= [s^2F_{M+1}^2(s)+2^2F_{M+2}(s)(F_{M+2}(s) - s F_{M+1}(s))] = s^2F_{M+1}^2(s)+2^2F_{M+2}(s)F_{M}(s) =A_0-2. \\ \nonumber
\end{align}
Hence Eq.~\ref{AutoFourier} and Eq.~\ref{Zeroq} are consistent with Eq.~\ref{HuffAutoCorrFinal}, as expected.  This check has revealed yet another interesting property of canonical Huffman arrays, namely that for such sequences the square of the sum equals the sum of the squares, bar the contributions from the $\pm1$ end points,
\begin{equation}\label{SumofSquares}
\sum_{i=2}^{N-1} [H^N_i(s)]^2 =[\sum_{i=2}^{N-1} H^N_i(s)]^2.
\end{equation}
Equation~\ref{SumofSquares} thereby shows that these canonically delta-correlated arrays correspond to an infinite family (indexed by $N$ and $s$) of Pythagorean n-tuples, whenever $s$ is an integer.  A related connection between Pell sequences (intrinsic here, when $s = 2$) and such Fermat-like multi-dimensional hypotenuse formulae has been noted in the literature \cite{PellPythagorasSumSquares2001}. These insights provide a simplified form of the power spectrum,
\begin{eqnarray}\label{PowerSpectrum}
|\mathcal{F}[H^N(s)]_{q}|^2 = [F_{M+2}(s) + F_{M}(s)]^2 + 2^2 \sin^2 (\pi q/N),
\end{eqnarray}
where we have substituted $2^2 \sin^2 (\pi q/N)$ for $2(1-\cos(2\pi q/N))$.
Equation~\ref{PowerSpectrum} shows that the mean of the Fourier magnitude is always less than $F_{M+2}(s) + F_{M}(s)$ and that the maximum can be reliably overestimated by setting $q=N$.  Using this along with a binomial expansion of the square root of Eq.~\ref{PowerSpectrum} leads to a succinct bound on the spectral-flatness $\mathcal{S}$ from Sec.~\ref{sec:deltaAuto},
\begin{equation}\label{PowerSpectrumBound}
    \mathcal{S}  <  2/[F_{M+2}(s) + F_{M}(s)]^2.
\end{equation}
Equation~\ref{PowerSpectrumBound} tends to zero in the large-$s$ limit, which can be seen from the fact that $|F_M(s)| \sim s^{M-1}$ for large $s$. Equation~\ref{PowerSpectrumBound} also tends to zero in the large-$M$ limit if $s \ge 1$, because the numerical value of any Fibonacci polynomial will become arbitrarily large if both of the following hold: (i) $M \gg 1$, (ii) $s\ge 1$.  This behavior of the spectral flatness $\mathcal{S}$ shows that the Fourier-spectra of Huffman arrays are asymptotically constant, just like the delta function.  Another way of expressing this important point is to consider the square root of Eq.~\ref{PowerSpectrum}. Using Eq.~\ref{ReconcileSumSquare} and the sine-square/cosine identify, this can be written as,
\begin{align}\label{FourierApprox}
|\mathcal{F}[H^N(s)]_q| = \sqrt{S_{N}^2 + 2^2 s_{q}^2} \approx |S_{N}| + (1-c_{q})/|S_{N}| -2 s_{q}^4/|S_{N}|^3 + 4 s_{q}^6/|S_{N}|^5 -10 s_{q}^8/|S_{N}|^7,
\end{align}
where $S_N$ is the sum over all elements of the Huffman array, $c_q$ is short for $\cos(2 \pi q/N)$ with $s_q = \sin(\pi q/N)$.  The leading term in Eq.~\ref{FourierApprox} is simply the magnitude of $S_N$ and this series arose from an expansion of the square root in Eq.~\ref{FourierApprox}, in terms of $1-c_{q}$.

To reach these conclusions we have evaluated the Fourier spectrum via indirect means, taking advantage of the Fourier correlation theorem to avoid the unnecessary complexities of exponential sums.  Due to the Golay skew symmetry of the Huffman arrays, it is however possible to evaluate the Fourier transform explicitly to yield both the magnitude and phase.

The discrete Fourier transform of the canonical Huffman arrays takes on a particularly simple form, if $H^N(s)$ is cyclically shifted such that the middle $\underline x$ term in $H^N(s)$ becomes the first element. With this rearrangement, the cyclically shifted sequence $H_c^N(s)$ almost coincides with the ordered Fibonacci polynomial sequence,
\begin{align}\label{HuffCycled}
H_c^N(s) = 2 s[&F_{M+1}(s)/2 - F_{M}(s)/s,F_{-M},F_{-(M-1)}(s), F_{-(M-2)}(s), \cdots,F_{-1}(s), -(2 s)^{-1},(2 s)^{-1}, F_{1}(s),\notag \\ &\cdots,F_{M-2}(s),F_{M-1}(s),F_{M}(s)],
\end{align}
with the exception of the first element and interior unit-magnitude elements.  Since spatial shifts of elements induce linear phase ramps in the spectral domain, this Huffman form remains canonically delta-correlated for \emph{periodic} auto-correlation.

Given the monotonically increasing indices, evaluating the Fourier transform of Eq.~\ref{HuffCycled} is almost the same task as that of the ordered Fibonacci polynomial sequence. Substitution of Eq.~\ref{HuffCycled} into Eq.~\ref{DiscreteFourier} gives a symmetrical form for the Fourier transform (using $M = (N-3)/2$),
\begin{align}\label{FourierHuffCycled}
\mathcal{F}[H_c^N(s)]_q = &s F_{(N-3)/2+1} - 2 F_{(N-3)/2} +2 s\sum_{n=1}^{(N-3)/2}F_{(N-1)/2-n}(e^{2 \pi i n q/N}+(-1)^n e^{-2 \pi i n q/N}) \notag \\ 
&+e^{-2 \pi i [(N+3)/2-1] q/N}-e^{2 \pi i [(N+3)/2-1] q/N},
\end{align}
where we have sifted out the $\pm 1$ Huffman elements from the Fourier sum.  Note the conjugated exponentials in Eq.~\ref{FourierHuffCycled} that alternately evaluate to sines or cosines.  Every element $\mathcal{F}[H^{N}(s)]_q$ can be written completely in terms of sinusoids, given known Chebyshev product-series for the Fibonacci polynomials \cite{LindCosines,SolutionFQCosines} or indeed trigonometric versions of the Binet form \cite{CahillFQCosines}. 

Using a computer algebra system, coupled with several of the aforementioned Fibonacci polynomial identities, after considerable algebraic reductions we were able to arrive at a compact closed-form expression for Eq.~\ref{FourierHuffCycled}, 
\begin{align}\label{FourierHuffFinal}
\mathcal{F}[H_c^N(s)]_q = -\frac{2 i \sin(2 \pi q/N)+s}{2 i \sin(2 \pi q/N)-s}[F_{2(N-3)/4}(s) +F_{2(N-3)/4+2}(s)+2 i (-1)^q \sin(\pi q/N)].
\end{align}
One can readily verify that Eq.~\ref{FourierHuffCycled} and Eq.~\ref{FourierHuffFinal} are equal using computer algebra constrained for real $s$ and $N = 4n+3$ for positive integer $n$.  Rather than pursuing a lengthy line-by-line derivation of Eq.~\ref{FourierHuffFinal}, it suffices to note that Eq.~\ref{FourierHuffCycled} can be expressed in closed-form since the Fibonacci polynomials can be written in terms of exponentials with phases that are linear in the summation index.  This arises from the well known connection between the Binet form and Chebyshev polynomials of the second kind \cite{CahillFQCosines}, such that
\begin{align}\label{BinetChebyshev}
F_{n}(s) &= i^{-n+1}\frac{e^{-i n \cos^{-1}(i s/2)}-e^{i n \cos^{-1}(i s/2)}}{e^{-i \cos^{-1}(i s/2)}-e^{i \cos^{-1}(i s/2)}}.
\end{align}
Hence the partial sum in Eq.~\ref{FourierHuffCycled} can be evaluated as variants of geometric series for the constituent terms, giving the periodic Fourier transform for  Fibonacci polynomial Huffman arrays in closed-form, albeit less succinct than that of Eq.~\ref{FourierHuffFinal}.

Equation~\ref{FourierHuffFinal} summarizes several of the key findings here, as the spectral flatness and canonical delta-correlation are  evident from the modulus.  Moreover, after cyclic permutation by $(N-1)/2+1$ steps, the inverse Fourier transform of Eq.~\ref{FourierHuffFinal} yields Huffman arrays, without reference to recursion or Binet forms, and hence may be viewed as a simpler representation of these canonical delta-correlated sequences. 

The scope of this paper has been to connect Hunt and Ackroyd's \cite{HuntAckroyd1980} Huffman sequences with the generalized Fibonacci sequence, to enable elementary and accessible proofs of the defining delta auto-correlation property.  It is anticipated that the closed-form Fourier transform expressions elucidated by this connection will lead to new insights concerning spectrally flat Huffman arrays.  As a concrete example, consider replacing the sines in the unit-amplitude phasor for the leading term in Eq.~\ref{FourierHuffCycled} with another periodic function that is odd in $q$, such as the tangent,  
\begin{align}\label{TangentHuffman}
\mathcal{F}[H_{t,cyclic}^N(s)]_q = 
-\frac{2 i \tan(2 \pi q/N)+s}{2 i \tan(2 \pi q/N)-s}[\left( \frac{2 +s}{2 -s} \right) ^{(N-1)/4}-\left(\frac{2 +s}{2 -s}\right)^{-(N-1)/4} +2 i (-1)^q \sin(\pi q/N)],
\end{align}
for scale parameter $|s| \ne 2$. The non-Fibonacci $q = 0$ term in Eq.~\ref{TangentHuffman} was arrived at by some educated guesses to ensure the canonical condition. Equation~\ref{TangentHuffman} has been checked numerically for $s$ ranging from -10 to +10 and $N = 4n+1$ from 3 to 45.  This was done by inverse Fourier transformation, the removal of redundant zeros in every second element, and subsequent cyclic permutation, the process of which gives rise to a scalable family of odd-length canonical Huffman arrays $H_t^N(s)$, where $N = 2n+1$. 

While Eq.~\ref{TangentHuffman} efficiently describes the Fourier spectra of cyclically permuted $H_t^N(s)$ (convenient for numerical evaluation), it is possible to state $H_t^N(s)$ directly as a family of sequences, for a particular choice of $s$.  To this end, setting the scale $s$ to unity yields the simple canonical form $H_{three}^N = [3, 3^{j - 1} 8,\cdots,\underline x ,\cdots,3^{-k - 1} 8, -1/3]$, with middle term $\underline x =3^{-(3-N)/2}-3^{(N-3)/2}$, where $j$ is the element number counted from the left and $k$ is that from the right of the array.  With only two non-zero off-peak auto-correlation elements equal to $-1$, these canonical Huffman sequences are not clearly defined by a recursion and do not exhibit skew symmetry.  Perhaps unsurprisingly, the zeros of the z-transform do not sit on golden-ratio circles with radii $X = \phi$ and $1/X = 1/\phi$ but instead we have $X = 3$ and $1/X = 1/3$.  For a given length $N$, the z-transforms from numerical inversion of Eq.~\ref{TangentHuffman} show that the complex arguments of the corresponding z-transform zeros sit at equi-angled phases (as required by Huffman's canonical conditions) while the radii $X(s), 1/X(s)$ vary with $s$, as expected. The diminishing fractions on the right hand side of this simple $H_{three}^N$ canonical sequence can in fact be truncated to derive yet another class of simple Huffman arrays. After rounding and truncating $H_{three}^N$, the following scalable canonical family of any length $N$ was inferred upon inspection,
\begin{align}\label{IntegerHuffman}
H_{int}^N(s) = [h_1, (s^2-1)s^n, \cdots, h_N].
\end{align}
Here, $h_1 = s$, $h_N = -s^{(N-2)}$ and the dots refer to the index $n$ running from $0$ to $N-3$.  The name of $H_{int}^N(s)$ refers to the all-integer valued elements that arise for integer choices of the real-valued parameter $s$. The canonical condition follows, for all lengths $N$ and finite scaling $s$, by appropriate summation of terms $H_{int, n}^N(s) H_{int, n+m}^N(s)$ for element indices $n$ and offsets $m$, including two more products to handle terms involving the $h_1$ and $h_N$ end elements. The auto-correlation peak is $1+s^{2N-2}$, with a pair of equal minima  $-s^{N-1}$ and all other terms zero.

 This work pursued elementary proofs of the canonical condition for Hunt and Ackroyd's \cite{HuntAckroyd1980} Fibonacci-based Huffman arrays.  The associated Fourier analysis was designed to prove and quantify the important spectral-flatness of the Fibonacci-based Huffman arrays but the simplicity of the approach has yielded unexpected insights, partly summarized by the families of canonical delta-correlated sequences $H_t^N(s)$, $H_{three}^N$ and $H_{int}^N(s)$.  Future work will further explore the solution space of canonical and quasi-Huffman sequences for greater flexibility in the design of such spectrally-flat delta-correlated arrays. 

\bibliographystyle{IEEEtran}

\bibliography{refs}

\begin{thebibliography}{10}
\providecommand{\url}[1]{#1}
\csname url@samestyle\endcsname
\providecommand{\newblock}{\relax}
\providecommand{\bibinfo}[2]{#2}
\providecommand{\BIBentrySTDinterwordspacing}{\spaceskip=0pt\relax}
\providecommand{\BIBentryALTinterwordstretchfactor}{4}
\providecommand{\BIBentryALTinterwordspacing}{\spaceskip=\fontdimen2\font plus
\BIBentryALTinterwordstretchfactor\fontdimen3\font minus
  \fontdimen4\font\relax}
\providecommand{\BIBforeignlanguage}[2]{{%
\expandafter\ifx\csname l@#1\endcsname\relax
\typeout{** WARNING: IEEEtran.bst: No hyphenation pattern has been}%
\typeout{** loaded for the language `#1'. Using the pattern for}%
\typeout{** the default language instead.}%
\else
\language=\csname l@#1\endcsname
\fi
#2}}
\providecommand{\BIBdecl}{\relax}
\BIBdecl

\bibitem{Singer1938}
J.~Singer, ``A theorem in finite projective geometry and some applications to
  number theory,'' \emph{Trans. Amer. Math. Soc.}, vol.~43, no.~3, pp.
  377--385, 1938.

\bibitem{McFarland2011}
G.~Oliveri, F.~Caramanica, C.~Fontanari, and A.~Massa, ``Rectangular thinned
  arrays based on {McFarland} difference sets,'' \emph{IEEE Trans. Antennas
  Propag.}, vol.~59, no.~5, pp. 1546--1552, 2011.

\bibitem{Luke1988}
H.~D. {Luke}, ``Sequences and arrays with perfect periodic correlation,''
  \emph{IEEE Trans. Aerosp. Electron.}, vol.~24, no.~3, pp. 287--294, 1988.

\bibitem{Huffman1962}
D.~{Huffman}, ``The generation of impulse-equivalent pulse trains,'' \emph{IRE
  Transactions on Information Theory}, vol.~8, no.~5, pp. 10--16, 1962.

\bibitem{Welch1974}
L.~{Welch}, ``Lower bounds on the maximum cross correlation of signals
  (corresp.),'' \emph{IEEE Trans. Inf.}, vol.~20, no.~3, pp. 397--399, 1974.

\bibitem{Sarwate1979}
D.~{Sarwate}, ``Bounds on crosscorrelation and autocorrelation of sequences
  (corresp.),'' \emph{IEEE Trans. Inf.}, vol.~25, no.~6, pp. 720--724, 1979.

\bibitem{Sarwate1999}
D.~V. Sarwate, ``Meeting the {W}elch bound with equality,'' in \emph{Sequences
  and Their Applications: Proceedings of SETA '98}, T.~H. C.~Ding and
  H.~Niederreiter, Eds.\hskip 1em plus 0.5em minus 0.4em\relax London:
  Springer-Verlag London Limited, 1999, pp. 79--102.

\bibitem{HuntAckroyd1980}
J.~{Hunt} and M.~{Ackroyd}, ``Some integer {H}uffman sequences (corresp.),''
  \emph{IEEE Trans. Inf.}, vol.~26, no.~1, pp. 105--107, 1980.

\bibitem{svalbe2019sharp}
I.~D. Svalbe, D.~M. Paganin, and T.~C. Petersen, ``Sharp computational images
  from diffuse beams: Factorization of the discrete delta function,''
  \emph{IEEE Trans. Comput. Imaging}, vol.~6, pp. 1258--1271, 2020.

\bibitem{Ackroyd1972}
M.~H. Ackroyd, ``Synthesis of efficient {H}uffman sequences,'' \emph{IEEE
  Trans. Aerosp. Electron. Syst.}, vol. AES-8, no.~1, pp. 2--8, 1972.

\bibitem{Ackroyd1975}
M.~{Ackroyd}, ``Huffman sequences with approximately uniform envelopes or
  cross-correlation functions (corresp.),'' \emph{IEEE Trans. Inf.}, vol.~23,
  no.~5, pp. 620--623, 1977.

\bibitem{Schroeder1970}
M.~{Schroeder}, ``Synthesis of low-peak-factor signals and binary sequences
  with low autocorrelation (corresp.),'' \emph{IEEE Trans. Inf}, vol.~16,
  no.~1, pp. 85--89, 1970.

\bibitem{Golay1972}
M.~J.~E. Golay, ``A class of finite binary sequences with alternate
  autocorrelation values equal to zero,'' \emph{IEEE Trans. Inf. Theory}, vol.
  IT-18, no.~3, pp. 449--450, 1972, correspondence.

\bibitem{Golay1975b}
M.~{Golay}, ``Hybrid low autocorrelation sequences (corresp.),'' \emph{IEEE
  Trans. Inf}, vol.~21, no.~4, pp. 460--462, 1975.

\bibitem{OjedaHuffman1994}
R.~Ojeda and E.~J. Tacconi, ``Huffman sequences with uniform time energy
  distribution,'' \emph{Signal Process.}, vol.~37, no.~1, pp. 141 -- 146, 1994.

\bibitem{White1977}
D.~J. White, J.~N. Hunt, and L.~A.~G. Dresel, ``Uniform {H}uffman sequences do
  not exist,'' \emph{B. Lond. Math. Soc.}, vol.~9, no.~2, pp. 193--198, 1977.

\bibitem{BorweinJedwab2004}
P.~{Borwein}, K.~.~S. {Choi}, and J.~{Jedwab}, ``Binary sequences with merit
  factor greater than 6.34,'' \emph{IEEE Trans. Inf.}, vol.~50, no.~12, pp.
  3234--3249, 2004.

\bibitem{Popovic1992}
B.~M. {Popovic}, ``Generalized chirp-like polyphase sequences with optimum
  correlation properties,'' \emph{IEEE Trans. Inf.}, vol.~38, no.~4, pp.
  1406--1409, 1992.

\bibitem{ChangGolombBarkerGener1996}
{N. Chang} and S.~W. {Golomb}, ``7200-phase generalized {B}arker sequences,''
  \emph{IEEE Trans. Inf.}, vol.~42, no.~4, pp. 1236--1238, 1996.

\bibitem{FriesePolyphBarker1996}
M.~{Friese}, ``Polyphase {B}arker sequences up to length 36,'' \emph{IEEE
  Trans. Inf.}, vol.~42, no.~4, pp. 1248--1250, 1996.

\bibitem{Golomb1982}
S.~{Golomb} and H.~{Taylor}, ``Two-dimensional synchronization patterns for
  minimum ambiguity,'' \emph{IEEE Trans. Inf.}, vol.~28, no.~4, pp. 600--604,
  1982.

\bibitem{Monk2004}
L.~Monk, D.~Tang, and D.~Brown, ``Identities for generalized {F}ibonacci
  numbers,'' \emph{Internat. J. Math. Ed. Sci. Tech.}, vol.~35, no.~3, pp.
  436--439, 2004.

\bibitem{Moharir1975}
P.~S. Moharir, ``Multilevel aperiodic {H}uffman sequences,'' \emph{Electron.
  Lett.}, vol.~11, no.~3, pp. 56--57, 1975.

\bibitem{PellPythagorasSumSquares2001}
R.~A. Beauregard and E.~R. Suryanarayan, ``Pythagorean boxes,'' \emph{Math.
  Mag.}, vol.~74, no.~3, pp. 222--227, 2001.

\bibitem{LindCosines}
D.~Lind, ``Problem {H}-64 proposal,'' \emph{Fibonacci Q.}, vol.~3, p. 116,
  1965.

\bibitem{SolutionFQCosines}
``Advanced problems and solutions,'' \emph{Fibonacci Q.}, vol.~5, pp. 74--75,
  1967.

\bibitem{CahillFQCosines}
N.~D. Cahill, J.~R. D’Errico, and J.~P. Spence, ``Complex factorizations of
  the {F}ibonacci and {L}ucas numbers,'' \emph{Fibonacci Q.}, vol.~41, pp.
  13--19, 2003.

\bibitem{Florez2018}
R.~Fl\'{o}rez, N.~McAnally, and A.~Mukherjee, ``Identities for the generalized
  {F}ibonacci polynomial,'' \emph{Integers}, vol. 18B, pp. 1--13, 2018.

\bibitem{JohnsonRuleFQ}
R.~Euler and J.~Sadek, ``Elementary problems and solutions,'' \emph{Fibonacci
  Q.}, vol.~42, no.~1, pp. 86--91, Feb. 2004.

\bibitem{Johnson}
\BIBentryALTinterwordspacing
R.~C. Johnson, \emph{Matrix methods for Fibonacci and related sequences},
  (accessed 7th May 2021). [Online]. Available:
  \url{https://maths.dur.ac.uk/~dma0rcj/PED/fib.pdf}
\BIBentrySTDinterwordspacing

\end{thebibliography}
\vskip -5\baselineskip
\begin{IEEEbiographynophoto}{Imants D. Svalbe}
completed a PhD in experimental nuclear physics at Melbourne University in 1979. His current work applies Mojette and Finite Radon transforms to design $nD$ geometric structures of signed integers that act as zero-sum projection ghosts in discrete tomography and to build large families of $nD$ integer arrays that have optimal correlation properties.
\end{IEEEbiographynophoto}
\vskip -5\baselineskip
\begin{IEEEbiographynophoto}{David M. Paganin}
received his PhD in optical physics from Melbourne University in 1999, and has been with Monash University since 2002. His research interests include x-ray optics, visible-light optics, electron diffraction, neutron optics and non-linear quantum fields. 
\end{IEEEbiographynophoto}
\vskip -5\baselineskip
\begin{IEEEbiographynophoto}{Timothy C. Petersen}
completed a PhD in the condensed matter physics of disordered carbon, at RMIT University in 2004.  He has performed experiments across a range of microscopy techniques to study disordered solids and develop new diffraction physics theories.
\end{IEEEbiographynophoto}

\section*{Supplementary Material}

In this Supplemental Material we give a formal proof that the aperiodic auto-correlation of the particular Huffman sequences derived by Hunt and Ackroyd \cite{HuntAckroyd1980} are canonical, using elementary identities for Fibonacci polynomials.   In other words, we show that such a Huffman sequence correlates as a discrete delta function: a peak at zero shift and exactly zero for all other shifts except the unavoidable end terms.  The proofs show how sums of products are zero for these off-peak shifts and the accompanying diagrams highlight the role of Golay skew-symmetry for canceling terms, the sign alternation of which is naturally provided by the Fibonacci property that  $f_{-j}  = ( - 1 )^{-j+1} f_{j}$ for every integer $j$. The first two short sections define the problem and the details of the proofs are contained in the longer Sec.~\ref{sec:Huffman proof}.  To supplement the introduction of the main text, the final section here (Sec.~\ref{sec:nonFibonacci}) provides a few canonical Huffman sequences that are unlike those of Hunt and Ackroyd \cite{HuntAckroyd1980} and which, though seemingly related, are not obviously based upon the Fibonacci sequence.  

\section{Delta Auto-Correlation} \label{sec: deltaAuto}
The $j^{th}$ element $A_j$ of an auto-correlation $A$  for any sequence is $S$ defined by
\begin{equation}\label{AutoCorr}
A_j = S\otimes S = \sum_{i} S_i S_{i+j},  
\end{equation}
 for shifts $|j|<N$ where the sum runs from elements $i = 1$ to the length of $S$, $i = N$, and we assume a convention that both $S_1$ and $S_N$ are finite. Shifts that take $i+j$ beyond $N$ contribute zero, which is equivalent to zero-padding the left of $S_1$ and right of $S_N$. The standard `aperiodic' auto-correlation definition in Eq.~\ref{AutoCorr} differs from the periodic version, whereby shifts $i+j$ beyond the length of $S$ correspond to cyclically permuted elements of $S$. The real-valued sequences we shall describe remain delta-correlated according to either definition.            

\section{Canonical Huffman Arrays} \label{sec:Huffman definiton}

For the largest misalignment of $S_i$ with $S_{i+j}$, the sum in Eq.~\ref{AutoCorr} constitutes a single element, hence $A_j$ must be non-zero for such shifts (there are a pair: one for the leftmost shift and another for the rightmost). Hence the closest possible representation of the delta function $\delta$ by $A$ must be of the form $[a,0,\cdots,0,A_0,0,\cdots,0,a]$ for the $2N-1$ auto-correlation elements. For integer valued sequences, we accordingly minimize the end elements such that $|a|= 1$. In the main text, delta-correlated sequences with corresponding $A$ that adhere to this ideal form are referred to as `canonical Huffman arrays' \cite{Huffman1962}.    

\section{Huffman arrays built from Fibonacci sequences are canonical } \label{sec:Huffman proof}

\begin{figure}
\centering
\includegraphics[ width=1.0\columnwidth]{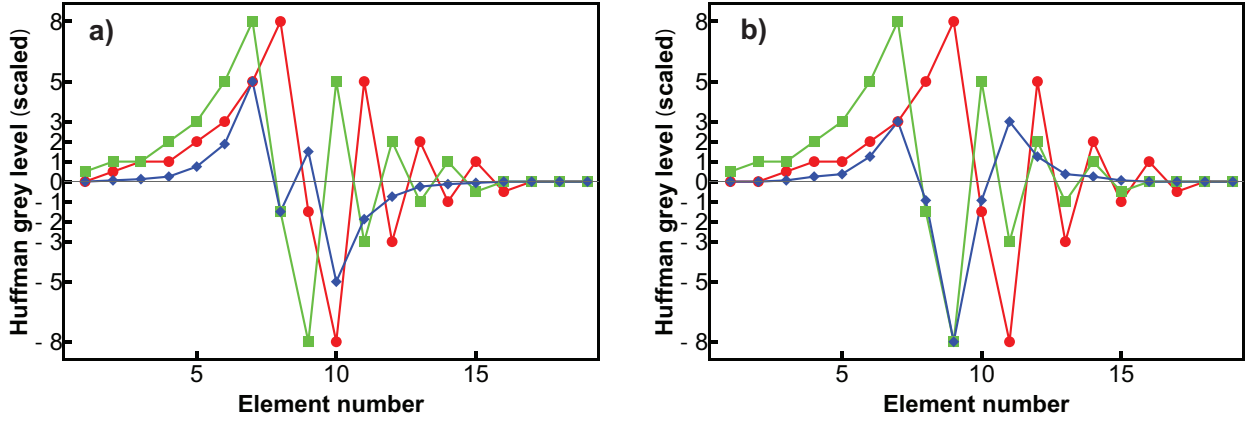} \caption{\textbf{Shifted products of Huffman arrays}. a) The canonical Huffman array $H^{15}$ was halved (green squares) then shifted by one element to the right (red circles) and the product of the two sequences was then divided by 8 (blue diamonds). Perfect anti-symmetry about the center is evident, implying a zero-sum for the auto-correlation at this shift. b) The same as for a), except that the halved $H^{15}$ Huffman sequence was instead shifted by 2 elements to the right. Unlike the graphs in a), it is difficult to see whether the sum of the product sequence cancels to zero or not.
}
\label{fig:HuffCancellation}
\end{figure}

The canonical delta-correlation condition is that all off peak elements in the auto-correlation $A$ must be zero, except for the end points. Golay \cite{Golay1972} has shown that all odd elements in $A$ are automatically zero, provided one generates $A$ from a seed sequence with `skew-symmetric' alternation of elements $S = [a,b,c,\cdots,c,-b,a]$ for lengths $4n+1$, where $n$ is an integer.  The Fibonacci based Huffman arrays here have length $4n+3$, for which the appropriate skew symmetry is instead $S = [a,b,c,\cdots,-c,b,-a]$.  Upon computing the auto-correlation of this seed sequence with integer elements and imposing the canonical condition, a set of quadratic Diophantine equations arises.  Minimizing the two end terms of the aperiodic auto-correlation $A$, such that $-a^2 = -1$, fixes $a$ to be $1$, upon which integer solutions arise when $b$ is even. The smallest choice $b = 2$ yields an odd-length sequence $H^N$ comprised of the negative and positive Fibonacci sequences, with an additional middle element $\underline{x}$,
%
\begin{align}\label{FibonacciHuffman}
H^N &= b[b^{-1},f_{1},f_{2},\cdots,f_{M},\underline{x}/b,f_{-M},\cdots,f_{-2},f_{-1},-b^{-1}] \\ \nonumber
&=[1,2f_{1},2f_{2},\cdots,2f_{M},\underline{x},2f_{-M},\cdots,2f_{-2},2f_{-1},-1]
\end{align}
%
where $M = (N-3)/2$, $\underline{x} = -f_{M-2}$ and $f_1, f_2, f_3, \cdots$ are elements of the Fibonacci sequence.  Integer solutions are only valid for even $M$, so we shall also impose this condition hereafter.  This restricts the length $N$ of the Huffman arrays to be of the form $N = 4n+3$ for integer $n \ge 1$.   
Fig.~\ref{fig:HuffCancellation}a demonstrates pictorially the cancellations that occur for odd shifts, which is clear from the product of $H^{15}/2$ (halved) with the same array shifted by 1 element (green squares and red circles, with the product shown as blue diamonds having been further scaled down by a multiplicative factor of 8). However the supposed cancellation for even shifts (2 elements here) is not clear at all in Fig.~\ref{fig:HuffCancellation}b, despite the fact that the summed product of the displaced graphs totals to zero. 

\begin{figure}[h!]
\centering
\includegraphics[ width=1.0\columnwidth]{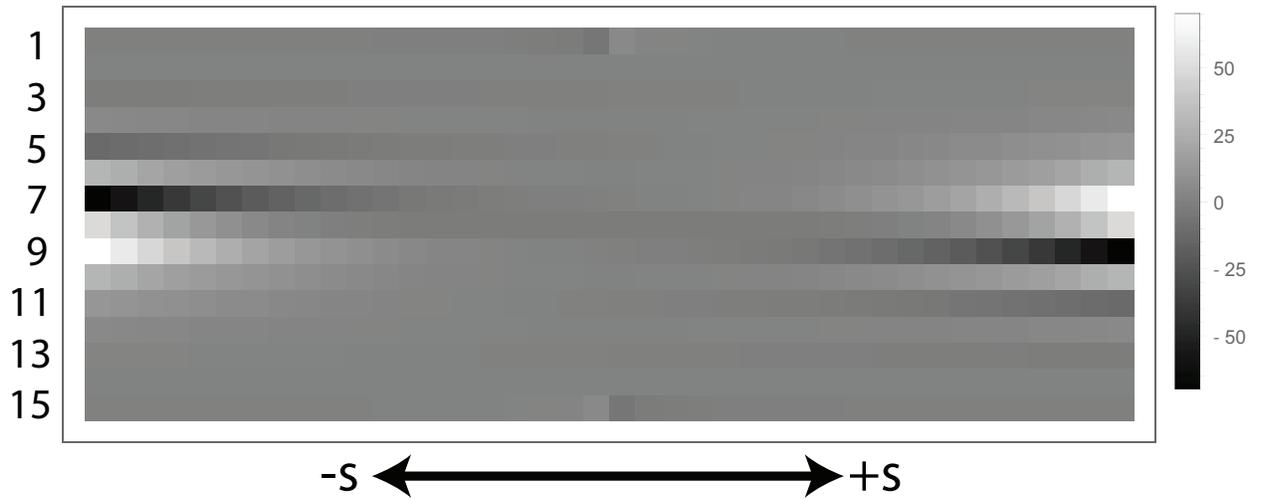} \caption{\textbf{Delta-correlated Huffman arrays of Fibonacci polynomials -- 15-element examples}. Discrete array elements of $H^{15}(s)/(2 s)$ correspond to vertical columns in the figure, where the odd element numbers have been labeled.  The division by $2 s$ highlights the detailed structure of the arrays and renders all elements (except the end points and middle point) as Fibonacci polynomial values.  All arrays are canonically delta-correlated and the horizontal row axis ranges from $s = -2$ to $s=2$, in steps of $0.1$, smoothly interpolating between the Fibonacci ($s=1$) and Pell ($s=2$) sequences.  Specifically, the $s=1$ case has the form $[1/2, 1, 1, 2, 3, 5, 8, -3/2, -8, 5, -3, 2, -1, 1, -1/2]$. The $s=2$ case corresponds to $[1/4, 1, 2, 5, 12, 29, 70, 99/2, -70, 29, -12, 5, -2, 1, -1/4]$. Although all canonical $H^N(s)$ arrays are well defined at $s=0$, the normalized $H^N(s)/(2 s)$ is not. Hence, the middle column for $s=0$ has been omitted.}
\label{fig:SmoothHuffman}
\end{figure}

The problem here is to show that Eq.~\ref{FibonacciHuffman} is canonically delta-correlated,
\begin{equation}\label{CanonicalAuto}
H^N\otimes H^N = [-1,0,\cdots,0,A_0,0,\cdots,0,-1],
\end{equation}
where $A_0 = \sum_{i=1}^{N} (H_i^N)^2$.  The generalization considered here is to prove Eq.~\ref{CanonicalAuto} for the cases where elements $f_1, f_2, \cdots ,f_j$ in Eq.~\ref{FibonacciHuffman} are replaced with Fibonacci polynomials $F_j(s)$, with the choice $b = 2 s$, and $s \in \mathbb{R}$. We refer to $s$ throughout as a `scaling parameter'.
For $s$ different from unity, the middle element generalizes to $\underline{x}= s F_{M+1}-2 F_{M}$.  The Fibonacci polynomials are defined by the recursion \cite{Florez2018},
\begin{equation}\label{FibPolynomial}
F_{r+2}(s) = s F_{r+1}(s)+F_{r}(s), 
\end{equation}
with $F_0(s) = 0$ and $F_1(s) = 1$.  For economy of notation, unless otherwise stated, we shall suppress the polynomial argument $s$ hereafter, such that $F_j(s) \equiv F_j$ and $H^N(s) \equiv H^N$.  Although the canonical Huffman arrays are designed to contain integer elements (for integer $s$), the delta-correlation property holds over continuous $s$. For example, $H^{11}(p/q)$ has entirely rational entries 
\begin{align}\label{RationalHuff}
H^{11}(m) = [1, 2 m, 2 m^2, 2 m (1 + m^2), 2 m^2 (2 + m^2), \\ \nonumber -3 m + m^3 + m^5,
-2 m^2 (2 + m^2), 2 m (1 + m^2), \\ \nonumber -2 m^2, 2 m, -1]
\end{align}
for $m=p/q$, where $p$ and $q$ are both non-zero integers.  As another example, $H^{11}(\phi)$ is irrational for all but the end points, with the interesting form %
\begin{align}\label{GoldenHuff}
H^{11}(\phi) = [1, 2 \phi, 2 \phi + 2, 6 \phi + 2, 10 \phi + 8, 
 2 (2 \phi + 2), \\ \nonumber -10 \phi - 8, 6 \phi + 2,- 2 \phi - 2, 
 2 \phi, -1],
\end{align}
where $\phi$ is the golden ratio.  Continuously scaled examples of Huffman arrays with $N=15$ are shown in Fig.~\ref{fig:SmoothHuffman}, over the range $s \in [-2, 2]$.

\begin{figure}
\centering
\includegraphics[ width=1.0\columnwidth]{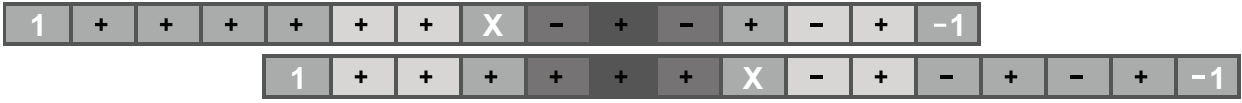}
\caption{\textbf{Auto-correlation of a skew-symmetric sequence for even shifts}. The $+$ and $-$ symbols indicate signs in the sequence $S$ from Eq.~\ref{AutoCorr}.  The X symbols represent the middle $\underline{x}$ element of both identical sequences.}
\label{fig:AutoCorrSmalld}
\end{figure}

Consider a specific auto-correlation element $A_d$ which can be deduced by shifting the bottom row of Fig.~\ref{fig:AutoCorrSmalld} relative to the top by $d$ positions ($d=4$ in this case), then multiplying aligned elements above and below before summing. 
On account of the skew-symmetric form of $H^N$, we need only examine the case where $d$ is even, 
\begin{align}\label{AutoCorrHuffman}
A_d = \sum_{i} H^N_d H^N_{i+d}.
\end{align}
First we consider correlation shifts $d$ of magnitude less than $M$ (the largest Fibonacci index). Next, the specific shift $d=M$ is treated separately.  The cases of larger shifts $d>M$ are dealt with in a subsequent sub-section (Sec.~\ref{sec:dGreaterthanM}). 

\subsection{Auto-correlation elements $A_d$ of $H^N(s)$ for $d<(N-3)/2$}

As $A_d = A_{-d}$, we need only consider even and positive values of $d$.  It is instructive to separate the product pairs that arise from Fibonacci terms with positive indices from those of the negative-index contributions, the two terms involving $\underline{x}$, as well as the terms arising from the unit-magnitude end points. Due to the imposed skew symmetry, the sum over positive-index pairs of terms is identical to a portion of the alternating sum, as portrayed by the lighter shaded elements in Fig.~\ref{fig:AutoCorrSmalld}.  Similarly, the alternating sum involving pairs of unmatched signs (comprising only one addend in the particular Fig.~\ref{fig:AutoCorrSmalld} example) is doubled by identical elements either side of the center of overlap. These paired elements are shaded darker in Fig.~\ref{fig:AutoCorrSmalld}, with the paired elements at the central overlap shaded darkest.  Lastly, the skew-symmetry ensures that the products involving $\pm$~1 are identical to each other, as are those involving the $\underline{x}$ elements. 

With these observations, element $A_d$ can then be summarized in terms of the Fibonacci polynomials as,
\begin{align}\label{AutoCorrSummary1}
A_{d}&= 8 s^2 [F_d/(2s)+(s F_{M+1}-2 F_{M}) F_{M-d+1}/(2s)\\ \nonumber&+(-1)^{r-1} F^2_{M-r+1}/2 
 +\sum_{i=1}^{M-d}F_i F_{i+d} \\ \nonumber &-\sum_{i=0}^{r-1}(-1)^i F_{M+i+2-d}F_{M-i}],
\end{align}
where $r = d/2$ and the third term corresponds to the darkest pair of elements in Fig.~\ref{fig:AutoCorrSmalld}.  Note that the first sum in Eq.~\ref{AutoCorrSummary1} is undefined for $d\geq M$. We will consider such larger shifts in the next two subsections. The pattern in Eq.~\ref{AutoCorrSummary1} for $d=0$ is consistent with $A_0 = \sum_{i=1}^N (H^N)^2$.

The first sum in Eq.~\ref{AutoCorrSummary1} comprises products of Fibonacci terms $F_{i}F_{i+d}$, which can be converted to a sum over squares using a recent generalization of  Catalan's identity that remains unchanged for Fibonacci polynomials \cite{Florez2018},
\begin{equation}\label{Catalan}
F_{n-r}F_{n+r} = F_{n}^2 - (-1)^{n-r} F_{r}^2.
\end{equation}
Shifting the dummy summation index by $r = d/2$ converts $F_{i}F_{i+d}$ into $F_{i-r}F_{i+r}$ such that
\begin{align}\label{CorrCatalan}
\sum_{i=1}^{M-d}F_i F_{i+d}&=\sum_{i=1+r}^{M-r}F_{i-r} F_{i+r} \\ \nonumber
&=\sum_{i=1+r}^{M-r} F_{i}^2 - (-1)^{i-r} F_{r}^2 \\ \nonumber
&=\sum_{i=1}^{M-r} F_{i}^2 - \sum_{i=1}^{r} F_{i}^2,
\end{align}
where an even-number of identical alternating terms were canceled to arrive at the final equality and the summation index was shifted again. The sums can be further simplified using the generalized Fibonacci sum of squares identity \cite{Monk2004},
\begin{equation}\label{SumSquares}
\sum_{i=1}^{N} F_i(s)^2 = F_N(s) F_{N+1}(s)/s.
\end{equation}
Hence, Eq.~\ref{CorrCatalan} becomes (suppressing the $s$ dependence again such that $F_j \equiv F_j(s)$),
\begin{equation}\label{CorrCatalan2}
\sum_{i=1}^{M-d}F_i F_{i+d}
=F_{M-r}F_{M-r+1}/s-F_{r}F_{r+1}/s,
\end{equation}
where the sum of squares identity was used twice.

The second summation in Eq.~\ref{AutoCorrSummary1} can be treated similarly,
\begin{align}\label{AutoCorrCatalan2}
&\sum_{i=0}^{r-1}(-1)^i F_{M+i+2-d}F_{M-i} \\ \nonumber =&\sum_{i=1-r}^{0}(-1)^{i-1+r} F_{M+1-r+i}F_{M+1-r-i}\\ \nonumber
=&\sum_{i=1-r}^{0}(-1)^{i-1+r} [F_{M+1-r}^2 - (-1)^{M+1-r-i} F_{i}^2].
\end{align}
The final sum in Eq.~\ref{AutoCorrCatalan2} can be written as
\begin{align}\label{AutoCorrCatalan3}
&\sum_{i=1-r}^{0}(-1)^{i-1+r} [F_{M+1-r}^2 - (-1)^{M+1-r-i} F_{i}^2] \\ \nonumber &= (1+(-1)^{r-1})/2 F_{M+1-r}^2 - \sum_{i=1}^{r-1} F_{i}^2,
\end{align}
where the factor of two accounts for the cancellation of all but one term, which occurs when $r$ is odd for the sum over $F_{M+1-r}$, and, by interchanging the upper and lower limits, we have used the fact that the negative-index Fibonacci terms have the same magnitude as those of the positive-index terms with $F_0(s)=0$. Using Eq.~\ref{SumSquares}, Eq.~\ref{AutoCorrCatalan2} becomes
\begin{align}\label{AutoCorrCatalan5}
&\sum_{i=0}^{r-1}(-1)^i F_{M+i+2-d}F_{M-i} 
\\ \nonumber &= (1+(-1)^{r-1})/2 F_{M+1-r}^2 - F_{r-1} F_{r}/s.
\end{align}
Pooling the equations above expresses the auto-correlation in closed-form,
\begin{align}\label{AutoCorrNoSums}
A_{d}&=8s^2[F_d/(2s)+(s F_{M+1}-2 F_{M}) F_{M-d+1}/(2s) \\ \nonumber
&+(-1)^{r-1} F^2_{M-r+1}/2 +F_{M-r}F_{M-r+1}/s-F_{r}F_{r+1}/s \\ \nonumber
&-(1+(-1)^{r-1}) F_{M+1-r}^2/2 + F_{r-1} F_{r}/s] \\ \nonumber
&=8s^2[F_d/(2s)+(s F_{M+1}-2 F_{M}) F_{M-d+1}/(2s) \\ \nonumber
&+F_{M-r}F_{M-r+1}/s-F_{r}F_{r+1}/s \\ \nonumber
&- F_{M+1-r}^2/2 + F_{r-1} F_{r}/s] \\ \nonumber
&=8s^2[F_{2r}/(2s)+(s F_{M+1}-2 F_{M})F_{M-2r+1}/(2s) \\ \nonumber
&-F_{M-r+1}^2/2 +F_{M-r}F_{M-r+1}/s-F_{r}^2],
\end{align} 
where the last equality was arrived at by using the Fibonacci polynomial recursion $F_{r+1} = s F_{r}+F_{r-1}$, such that $F_r(F_{r+1}-F_{r-1})/s = F_{r}^2$. Also, $d$ is now notated as $2r$ hereafter.

None of the terms in Eq.~\ref{AutoCorrNoSums} cancel in pairs, for arbitrary indices, hence further transformations are required at this point. The `Johnson rule' for product pairs \cite{JohnsonRuleFQ}, from which many other second-order Fibonacci number relations can be derived \cite{Johnson}, is of great utility here:  
\begin{equation}\label{Johnson}
f_{a}f_{b} - f_{c}f_{d} = (-1)^t(f_{a-t}f_{b-t} - f_{c-t}f_{d-t}),
\end{equation}
where $a,b,c,d,t$ are any integers and the constraint $a+b = c+d$ must be satisfied.  This rule remains unchanged for Fibonacci polynomials \cite{JohnsonRuleFQ} (whereby all $f$ are replaced by $F$ in Eq.~\ref{Johnson}).  There are two pairs of consecutive terms in Eq.~\ref{AutoCorrNoSums} with indices involving both $M$ and $r$ that match this form and meet this requirement. Before employing the Johnson rule, note that the $t$ parameter can be used to match any of the indices and thereby set one of the constituent pairs to zero (since $f_0 = 0$).  With this in mind, choosing $t = M-2r+1$ removes two terms, so that Eq.~\ref{AutoCorrNoSums} becomes
\begin{align}\label{AutoCorrNoSums2}
A_{2r}
&=8s^2[F_{2r}/(2s)
+F_{r}^2/2 -F_{r-1}F_{r}/s-F_{r}^2] \\ \nonumber
&=8s^2[F_{2r}/(2s)
 -F_{r-1}F_{r}/s-F_{r}^2/2], 
\end{align} 
where $t$ is odd since $M$ is even. To complete the proof that $A_d = 0$ for $0 < d < M$, we will use the identity (which can be verified by substituting the Binet forms on both sides \cite{Florez2018}),
\begin{equation}\label{FlorezProp1}
F_{m+n+1} = F_{m+1}F_{n+1}+ F_{m} F_{n}.
\end{equation}
Setting $m=r$ and $n = r-1$ in Eq.~\ref{FlorezProp1}, we now have,
%
\begin{align}\label{AutoCorrNoSums3}
A_{2r}
&=8s^2[(F_{r+1}F_{r}+ F_{r}F_{r-1})/(2s)
 -F_{r-1}F_{r}/s-F_{r}^2/2] \\ \nonumber 
&=8s^2[F_{r+1}F_{r}/(2s)
 -F_{r-1}F_{r}/(2s)-F_{r}^2/2].
\end{align} 
%
The Fibonacci polynomial identity $F_{r-1} = s F_{r}+F_{r-1}$ dictates that $(F_{r+1}-F_{r-1})/s=F_{r}$.  Using this gives,
\begin{equation}\label{AutoCorrNoSums4}
A_{2r}=8s^2[F_{r}F_{r}/2-F_{r}^2/2]=0.
\end{equation} 
Hence, for cases where $d<M$, we have proven that Fibonacci polynomial based Huffman sequences are canonical for all even $M$ and any scaling parameter $s$.

\subsection{Auto-correlation element $A_d$ of $H^N(s)$ for $d=(N-3)/2$}

As Fig.~\ref{fig:AutoCorrSmalld} implies, the sum over strictly positive-index Fibonacci pairs does not occur at shift $d=M$, so the first summation in Eq.~\ref{AutoCorrSummary1} must be omitted in this specific case.  Upon substituting $d=M$, the absence of this term renders Eq.~\ref{AutoCorrNoSums} as the simpler form,
%
\begin{equation}\label{AutoCorrShiftisM}
A_{M}=8s^2[(s F_{M+1}-F_{M})/(2 s) + F_{M/2-1}F_{M/2}/s -F_{M/2+1}^2/2].
\end{equation}
%
To show this is zero, it is convenient to match the indices in which $M$ is halved, with those where $M$ is not halved. Writing $t=M/2$, this can be facilitated using both Catalan's identity Eq.~\ref{Catalan} to give $F_{2t+1} = F_{t}^2+F_{t+1}^2$ as well as Eq.~\ref{FlorezProp1} in the form $F_{2t}=F_{t}(F_{t-1}+F_{t+1})$.  Using these expressions, Eq.~\ref{AutoCorrShiftisM} becomes    
\begin{align}\label{AutoCorrShiftisM}
A_{2t}&=8s^2[(F_{t}^2+F_{t+1}^2)/2-F_{t}(F_{t-1}+F_{t+1})/(2 s) \\ \nonumber
&+ F_{t-1}F_{t}/s -F_{t+1}^2/2] \\ \nonumber
&=8s^2[F_{t}^2/2 - F_{t}(F_{t+1} + F_{t-1})/(2 s)] \\ \nonumber
&=0,
\end{align} 
where the last line is from the Fibonacci polynomial recursion on the elements in round brackets.

\subsection{Auto-correlation elements $A_d$ of $H^N(s)$ for $d > (N-3)/2$}
\label{sec:dGreaterthanM}
Consider now the remaining case for larger shifts $d > M$, which do not involve the middle $\underline{x}$ elements.  
\begin{figure}[h!]
\centering
\includegraphics[ width=1.0\columnwidth]{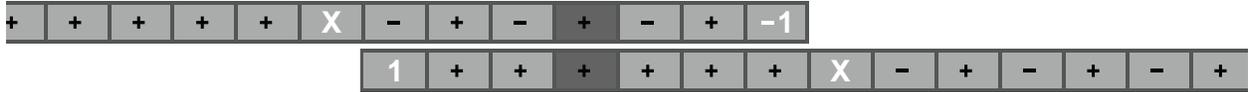}
\caption{\textbf{Auto-correlation $A_d$ of a skew-symmetric sequence for even shifts with $d> M$}. The $+$ and $-$ symbols indicate signs in the sequence $S$ from Eq.~\ref{AutoCorr}.  The X symbols represent the middle $\underline{x}$ element of both identical sequences.  Some unpaired elements have been cropped out of the field of view, on the left and right of the schematic.}
\label{fig:AutoCorrLarged}
\end{figure}
Upon inspection of Fig.~\ref{fig:AutoCorrLarged}, we can infer that the auto-correlation $A_d$ comprises two identical sums over pairs of positive- and negative-index Fibonacci polynomials, either side of the center of overlap, as well as two identical terms involving $\pm~1$. The darkest shaded pair of (identical) elements in Fig.~\ref{fig:AutoCorrLarged} gives rise to a solitary squared term in $A_d$.  The sum of all these terms is given in Eq.~\ref{AutoCorrLarged} (with $d = 2r$),
\begin{align}\label{AutoCorrLarged}
A_{2r}=8 s^2 [-F_{2(M-r+1)}/(2s)+(-1)^{r}F_{M-r+1}^2/2 \\ \nonumber -\sum_{i=1}^{M-r}(-1)^{i}F_i F_{2(M-r+1)-i}].
\end{align} 
We can simplify the summation in Eq.~\ref{AutoCorrLarged} using similar procedures to those employed earlier. To this end, Catalan's identity Eq.~\ref{Catalan} can be used if the index in the summation is first shifted, 
%
\begin{align}\label{AutoLargedSum}
&\sum_{i=1}^{M-r}(-1)^{i}_i F_{2(M-r+1)-i} \\ \nonumber
= &\sum_{i=-(M-r)}^{-1}(-1)^{i} F_{(M-r+1)+i} F_{(M-r+1)-i},
\end{align} 
%
in which the addends can be written, 
%
\begin{align}\label{AutoLargedAddends}
&F_{(M-r+1)+i}F_{(M-r+1)-i} \\ \nonumber
&= F_{M-r+1}^2 - (-1)^{(M-r+1)-i} F_{i}^2.
\end{align}
%
Inserting Eq.~\ref{AutoLargedAddends} into Eq.~\ref{AutoLargedSum} and also shifting the index back by $M-r+1$ yields,
\begin{align}\label{AutoLargedSumII}
&\sum_{i=1}^{M-r}(-1)^{i}F_i F_{2(M-r+1)-i} \\ \nonumber
&=\sum_{i=1}^{M-r}(-1)^{i}(F_{M-r+1}^2 - (-1)^{-i} F_{i-(M-r+1)}^2) \\ \nonumber
&=-(1-(-1)^{r})/2 F_{M-r+1}^2 -\sum_{i=1}^{M-r}F_{i-(M-r+1)}^2 \\ \nonumber
&=-(1-(-1)^{r})/2 F_{M-r+1}^2 -\sum_{i=-(M-r)}^{-1}F_{i}^2,
\end{align} 
where $(-1)^{2(M-r+1)} =1$ and $(-1)^{M-r} = (-1)^r$ were used (for the latter, $M$ is even) and all alternating signed terms of constant magnitude cancel when $r$ is even, leaving one final term (inducing the division by 2).

Since the magnitude of negative indices matches that of positive indices, the limits of the sum in Eq.~\ref{AutoLargedSumII} can be interchanged to reduce this to a single term using the generalized sum over squares identity Eq.~\ref{SumSquares},
\begin{align}\label{AutoLargedSumIII}
&\sum_{i=1}^{M-r}(-1)^{i}F_i F_{2(M-r+1)-i} \\ \nonumber
&= -(1-(-1)^{r})/2 F_{M-r+1}^2-F_{M-r}F_{M-r+1}/s.
\end{align} 
Inserting Eq.~\ref{AutoLargedSumIII} into Eq.~\ref{AutoCorrLarged} yields the auto-correlation $A_{2r}$ at even shift $d$ as a closed-form expression comprising three terms,
%
\begin{align}\label{AutoCorrLargedClosed}
&A_{2r}=8 s^2 [-F_{2(M-r+1)}/(2s) +(-1)^{r}F_{M-r+1}^2/2 \\ \nonumber
&+(1-(-1)^{r})/2 F_{M-r+1}^2 + F_{M-r}F_{M-r+1}/s \\ \nonumber
&= - 8 s^2 [F_{2(M-r+1)}/(2s)
- F_{M-r}F_{M-r+1}/s-F_{M-r+1}^2/2],
\end{align} 
%
where we have factored out $-1$ to reveal similarity with Eq.~\ref{AutoCorrNoSums}.  By comparison, the terms between square brackets are identical if the Fibonacci index $r$ in Eq.~\ref{AutoCorrNoSums} is substituted with the index $M-r+1$.  As such Eq.~\ref{AutoCorrLargedClosed} is automatically zero (i.e.~$A_d=0$, for $d > M$). This completes the proof that Huffman arrays, constituted by skew-symmetric Fibonacci polynomial sequences, are canonical for all auto-correlation shifts.

\section{Other canonical sequences not based on the Fibonacci recursion} 
\label{sec:nonFibonacci}
There are many canonical Huffman sequences which can be found by direct methods that do not follow the recursion for Fibonacci polynomials, some comprising integer elements.  Though outside the scope of this work, we provide three such examples for the interested reader, which are different from those discussed in the main manuscript.  The following length 11 non-Fibonacci canonical Huffman array $H_{non}^{11}$ exhibits all the same correlation metrics as the Fibonacci based $H^{11}(s=1)$,
\begin{align}\label{H11nonFib}
H_{non}^{11} = [1, 1, 3, 4, 2, 6, -7, -1, 2, 1, -1].
\end{align} 
Likewise, the z-transform zeros sit on circles in the Argand plane with the same $\phi$ and $1/\phi$ radii as that of $H^{11}(s=1)$, where $\phi$ is the golden ratio.  

There exist canonical arrays of length $4n+1$ and, as remarked by Hunt and Ackroyd \cite{HuntAckroyd1980}, no entirely-integer solutions exist when the end elements have opposite signs $\pm1$.  However, for arrays with sufficiently large elements, length $4n+1$ canonical arrays which asymptotically approximate integer elements, for example this almost-integer non-Fibonacci based Huffman array $H_{non}^{9}$,
\begin{align}\label{H9}
&H_{non}^{9} = [1, 200, 100 (200 - 2 \sqrt{10002}),  \\ \nonumber
&100 (-2 - 400 \sqrt{10002}), -4000000 \sqrt{10002},  \\ \nonumber
&100 (-2 + 400 \sqrt{10002}), 100 (-200 - 2 \sqrt{10002}), 200, -1]
\end{align} 
Also, it is possible to find canonical length $4n+1$ arrays comprising Gaussian integer elements, such as $[1, 2i, -2 - 2i, 4 - 2i, 4i, -4 - 2i, 2 - 2i, 2i, -1]$

Lastly, if we consider the Huffman extension suggested by Hunt and Ackroyd \cite{HuntAckroyd1980} whereby both end elements are matched to be $+1$, then integer $4n+1$ canonical arrays are again possible, such as,  
\begin{align}\label{H13}
&H_{non}^{13} = [1, 4, 8, 14, 24, 20, -14, -20, 24, -14, 8, -4, 1].
\end{align} 

\end{document}